\documentclass[a4paper]{article}
\usepackage{amssymb}
\usepackage[polish,english]{babel}

\date{\empty}
\newcommand{\Zn}{\mathbb{Z}}
\newcommand{\Rn}{\mathbb{R}}
\newcommand{\Cn}{\mathbb{C}}

\newcommand{\ord}{{\rm ord}\,}
\newcommand{\scalar}[2]{{\langle #1, #2\rangle}}

\newtheorem{Theorem}{Theorem}
\newtheorem{Lemma}[Theorem]{Lemma}
\newtheorem{Corollary}[Theorem]{Corollary}

\begin{document}

\title{Kouchnirenko type formulas for local invariants of plane analytic curves}
\author{\foreignlanguage{polish}{Janusz Gwo"zdziewicz}
}
\maketitle

\begin{abstract}
Let $f(x,y)=0$, $g(x,y)=0$ be equations of plane analytic curves defined in the 
neighborhood of the origin and let $\pi:M\to(\Cn^2,0)$ be a local toric modification.
We give a formula which connects a number of double points hidden at zero $\delta_0(f)$
with a sum $\sum_p \delta_p(\tilde f)$ which runs over all intersection points of the 
proper preimage of $f=0$ with the exceptional divisor. We give also similar formulas 
for the Milnor number $\mu_0(f)$ and the intersection multiplicity $(f,g)_0$.
Presented formulas generalize Kouchnirenko and Bernstein theorems and classical Noether
formula for the intersection multiplicity after blow-up. 
\end{abstract}

\section{Kouchnirenko and Bernstein theorems}

This paper is an english translation of \cite{Gw}.

Let $f\in\Cn\{x,y\}$ be a convergent power series such that
$f(x,0)$ and $f(0,y)$ do not vanish. Such a series is called convenient.
For any convenient series $f(x,y)=\sum_{ij}a_{ij}x^iy^j$ 
the Newton diagram $\Delta_f$ is the convex closure of the set
$$\bigcup_{a_{ij}\neq0}\{\,(i,j)+ \Rn_{+}^2\,\} .
$$
Since  $f$ is convenient,
$\Delta_f$ has joint points with both axis.
A union of compact edges of $\Delta_f$ is called
the Newton polygon of~$f$ and is denoted ${\cal N}_f$.

Introduce further notations. For the Newton diagram $\Delta$ whose Newton polygon ${\cal N}$
touches axes at points $(a,0)$ and $(0,b)$ denote:

\begin{itemize}
\item $P(\Delta) =\mbox{Area}\,(\Rn_{+}^2\setminus\Delta)$
\item $\mu(\Delta)=2P(\Delta)-a-b+1$
\item $r(\Delta)={}$ number of lattice points on ${\cal N}$ minus 1
\item $\delta(\Delta)=(\mu(\Delta)+r(\Delta)-1)/2$
\end{itemize}

If $\Delta_1$, $\Delta_2$ are two Newton diagrams then their mixed Minkowski volume 
is the quantity $[\Delta_1,\Delta_2]=P(\Delta_1+\Delta_2)-P(\Delta_1)-P(\Delta_2)$.

In 1970´s mathematicians from Arnold´s seminar
gave many formulas for invariants of singularities in terms of Newton diagrams.
We quote some of their results (see.~\cite{Ko}, \cite{Kh}).

\begin{Theorem}\label{T1}
Let $f,g\in\Cn\{x,y\}$ be convenient convergent power series. Then
\begin{equation}\label{T1.1} (f,g)_0\geq [\Delta_f,\Delta_g]
\end{equation}
\begin{equation}\label{T1.2} \delta_0(f)\geq \delta(\Delta_f)
\end{equation}
\begin{equation}\label{T1.3} \mu_0(f)\geq \mu(\Delta_f)
\end{equation}
\begin{equation}\label{T1.4} r_0(f)\leq r(\Delta_f)
\end{equation}
If coefficients of $f(x,y)$, $g(x,y)$, with indices from the sets
${\cal N}_f$, ${\cal N}_g$ respectively, are sufficiently general 
then (1)--(4) become equalities.
\end{Theorem}

In Theorem~\ref{T1}: 
$(f,g)_0$~is an intersection number of~$f$ and~$g$ at the origin,
$\delta_0(f)$ is the number of double points of a curve~$f=0$ hidden at zero,
$\mu_0(f)$~is the Milnor number of $f$ at zero and
$r_0(f)$~is a number of branches of $f=0$ at 0. 
``Sufficiently general'' in the second part of the theorem 
means that coefficients of power series  $f(x,y)$, $g(x,y)$ with indices
from the sets ${\cal N}_f$, ${\cal N}_g$ respectively
(there is a finite number of such coefficients) belong to some dense constructible set.
The equations of this set form so called nondegeneracy conditions (see \cite{Kh}).

\section{Noether theorem}
Let $f(x,y)=0$, $g(x,y)=0$ be equations of analytic curves defined in the neighborhood 
of 0 of the complex plane. 
Let $\sigma:M\to\Cn^2$ be a blowing-up of $\Cn^2$ at 0.
There are classical formulas which connect invariants of singularities of these curves 
with invariants of singularities of their proper preimages.

\begin{Theorem}\label{T2}
Assume that curves $f=0$, $g=0$ do not have common components. If $\tilde f=0$, $\tilde g=0$
are local equations of their proper preimages under blowing-up $\sigma$ then
\begin{equation}\label{T2.1}
(f,g)_0 = (\ord f)(\ord g) + \sum_{p\in\sigma^{-1}(0)} (\tilde f, \tilde g)_p
\end{equation}
If a curve $f=0$ does not have multiple components then
\begin{equation}\label{T2.2}
\delta_0(f)=(\ord f)(\ord f-1)/2 + \sum_{p\in\sigma^{-1}(0)} \delta_p(\tilde f)
\end{equation}
\begin{equation}\label{T2.3}
\mu_0(f)-1=(\ord f)(\ord f-1) + \sum_{p\in\sigma^{-1}(0)\cap\{\tilde f=0\}} (\mu_p(\tilde f)-1)
\end{equation}
\end{Theorem}

The classical Noether theorems are formulas~(\ref{T2.1}) for the intersection multiplicity
(\cite{B}, Theorem~13) and~(\ref{T2.2}) for the number of double points~$\delta_0(f)$.
A formula~(\ref{T2.3}) for the Milnor number follows from equality $2\delta_0(f)=\mu_0(f)+r_0(f)-1$.

\section{Local toric modifications}
The aim of this paper is to give formulas which generalize these from Theorem~\ref{T2} to
the case of an arbitrary local toric modification~$\pi:M\to\Cn^2$. 
When~$\pi$ will be a blow-up they will reduce to (\ref{T2.1})--(\ref{T2.3}) 
and for a toric modification with a sufficiently subtle fan they will give (\ref{T1.1})--(\ref{T1.3}) 
from Theorem~\ref{T1}.

\subsection{Fans}
By a simple cone $\sigma\subset\Rn^2$ we mean the set 
$$
\sigma=\sigma[\vec\xi,\vec\nu]=\{\,\alpha\vec\xi+\beta\vec\nu: \alpha\geq0, \beta\geq0\,\}
$$
where vectors $\vec\xi$, $\vec\nu\in\Rn^2$ have integer coordinates 
and form a base of lattice~$\Zn^2$, i.e. $\det(\vec\xi, \vec\nu)=\pm1$.

By a fan we mean a finite set of simple cones such that their union is the first 
quadrant~$\Rn_{+}^2$ and such that they intersect at most along edges. 
For every fan~$\cal W$ consisting of~$n$ cones one can enumerate 
counter clockwise the shortest lattice vectors from its rays. 
We get a sequence $\vec\xi_0$, $\vec\xi_1$, \dots, $\vec\xi_n$, 
where $\vec\xi_0=(1,0)$, $\vec\xi_n=(0,1)$
and $\det(\vec\xi_{i-1}, \vec\xi_i)=1$ for $1\leq i\leq n$.
We will say that $\cal W$ is spanned by $\vec\xi_0$, $\vec\xi_1$, \dots, $\vec\xi_n$.

\subsection{Local toric modifications}
With every simple cone $\sigma=\sigma[\vec\xi,\vec\nu]$, where 
$\vec\xi=(\xi_1,\xi_2)$, $\vec\nu=(\nu_1,\nu_2)$, $\det(\vec\xi, \vec\nu)=1$,
we associate a mapping $\varphi_\sigma:\Cn^2\to\Cn^2$ given in coordinates by
$$\varphi_\sigma:\left\{
\begin{array}{rcl}
x &=& u^{\xi_1}v^{\nu_1} \\
y &=& u^{\xi_2}v^{\nu_2} \\
\end{array}
\right.
$$
Let $\cal W$ be a fan consisting of $n$ cones
\begin{Theorem}
There exist a smooth analytic manifold $M$ and a proper analytic mapping
$\pi:M\to\Cn^2$ such that:
\begin{description}
\item[(i)] $\pi$ is an isomorphism from $M\setminus \pi^{-1}(0)$ to $\Cn^2\setminus \{0\}$,

\item[(ii)] the manifold $M$ is covered by $n$ charts associated with cones $\sigma_i$ ($i=1\dots n$)
            and in local coordinates of $i$-th chart the mapping~$\pi$ is given by formula $\pi=\varphi_{\sigma_i}$.
\end{description}
\end{Theorem}
We call $\pi:M\to\Cn^2$ a local toric modification associated with a fan $\cal W$.

\subsection{Thickened Newton diagrams}
For a Newton diagram~$\Delta$ and $\vec\xi\in \Rn_{+}^2$ we define 
a support function
$$
   l(\Delta,\vec\xi) = \inf_{p\in\Delta} \scalar{p}{\vec\xi}
$$

For a fan $\cal W$ spanned by vectors
$\vec\xi_0$, $\vec\xi_1$, \dots, $\vec\xi_n$ and
a Newton diagram~$\Delta$
we define $\tilde\Delta$ as an intersection of $n+1$ half-planes
$$
 \tilde\Delta = \bigcap_{i=0}^n
    \{\, p\in\Rn_{+}^2: \scalar{p}{\vec\xi_i} \geq l(\Delta,\vec\xi_i) \,\}
$$
and call this set a thickened Newton diagram relative to~$\cal W$ .
It follows directly from definition that
$l(\tilde\Delta,\vec\xi_i)=l(\Delta,\vec\xi_i)$ for $i=0$, \dots, $n$.

\section{Generalized Kouchnirenko theorem}

\begin{Theorem}\label{T3}
Let $\pi:M\to\Cn^2$ be a local toric modification associated with a fan~$\cal W$. 
If $f$, $g\in\Cn\{x,y\}$ are convenient power series 
and $\tilde f$, $\tilde g$ are their proper preimages then
\begin{equation}\label{T3.1}
  (f,g)_0 =
    [\tilde\Delta_f,\tilde\Delta_g] +\sum_{p\in\pi^{-1}(0)} (\tilde f,\tilde g)_p
\end{equation}
\begin{equation}\label{T3.2}
  \delta_0(f) = \delta(\tilde\Delta_f) + \sum_{p\in\pi^{-1}(0)} \delta_p(\tilde f)
\end{equation}
\begin{equation}\label{T3.3}
  \mu_0(f) = \mu(\tilde\Delta_f)+ r(\tilde\Delta_f) +
       \sum_{p\in\pi^{-1}(0)\cap\{\tilde f=0\}} (\mu_p(\tilde f)-1)
\end{equation}
\end{Theorem}

To avoid considering special cases in the statement of the theorem 
we adopt usual conventions about adding $+\infty$.

\vspace{1ex}
\textbf{Example~1.} The simplest fan $\cal W$ has only one cone~$\sigma$ 
which is the first quadrant. It is spanned by vectors
$\vec\xi_0=(1,0)$, $\vec\xi_1=(0,1)$. A mapping $\varphi_\sigma$ is given by
formula $(x,y)=(u,v)$ so the local toric modification associated with~$\cal W$
is an identity $\Cn^2\to\Cn^2$.

If $\Delta$ is a Newton diagram of a convenient power series then $\tilde\Delta=\Rn_{+}^2$. Hence the
invariants of a thickened diagram are:
$\mu(\tilde\Delta)=1$, $r(\tilde\Delta)=0$, $\delta(\tilde\Delta)=0$, $[\tilde\Delta,\tilde\Delta]=0$
and Theorem~\ref{T3} is trivially satisfied.

\vspace{1ex}
\textbf{Example~2.}  Consider a fan $\cal W$ consisting of two cones.
It is easy to check that~$\cal W$ is spanned by: 
$\vec\xi_0=(1,0)$, $\vec\xi_1=(1,1)$, $\vec\xi_2=(0,1)$.
The first cone $\sigma_1$ is generated by vectors $\vec\xi_0=(1,0)$, $\vec\xi_1=(1,1)$ 
and the mapping $\varphi_{\sigma_1}$ is given by a formula $(x,y)=(uv,v)$.
A mapping~$\varphi_{\sigma_2}$ associated with the second cone is given by $(x,y)=(u,uv)$.
Hence the toric modification $\pi:M\to\Cn^2$ associated with~$\cal W$ 
is a blowing-up of~$\Cn^2$ at zero.

If $\Delta$ is a Newton diagram of a convenient power series~$f$ of order~$d$
then the Newton polygon of~$\tilde\Delta$ is a segment with endpoints $(d,0)$ and $(0,d)$.
Hence
$P(\tilde\Delta)=d^2/2$,\quad
$\mu(\tilde\Delta)=2P(\tilde\Delta)-d-d+1=(d-1)^2$,\quad
$r(\tilde\Delta)=d$,\quad
$\delta(\tilde\Delta)=(\mu(\tilde\Delta)+r(\tilde\Delta)-1)/2=d(d-1)/2$.
If $f$, $g$ are convenient power series then $[\tilde\Delta_f,\tilde\Delta_g]=(\ord f)(\ord g)$. 
Substituting these quantities to equations~(\ref{T3.1})--(\ref{T3.3}) 
from Theorem~\ref{T3} we see that they reduce to~(\ref{T2.1})--(\ref{T2.3}) from Theorem~\ref{T2}.

\vspace{2ex}
We checked that Theorem~\ref{T2} is a special case of Theorem~\ref{T3}.
Likewise Theorem~\ref{T1}. If $f$, $g$ are convenient power series, then
there exists such a fan $\cal W$, that among its  spanning vectors are vectors 
orthogonal to all segments of Newton polygons~${\cal N}_f$ and~${\cal N}_g$. 
Then~$\tilde\Delta_f=\Delta_f$ and~$\tilde\Delta_g=\Delta_g$.
Inequalities in Theorem~\ref{T1} follow from adding extra terms on the right-hand side of~(\ref{T3.1})--(\ref{T3.3}).
If $f$ and $g$ satisfy nondegeneracy conditions then their proper preimages 
define smooth curves which do not have joint points on exceptional divisor~$\pi^{-1}(0)$
(see~\cite{Kh}), the sums on the right-hand side of~(\ref{T3.1}) and~(\ref{T3.2}) are 0 and a 
sum on the right-hand side of~(\ref{T3.3}) equals $-r_0(f)=-r(\Delta_f)$. Hence also
Theorem~\ref{T1} is a special case of Theorem~\ref{T3}.

\subsection{A decomposition of a toric modification to blowing-ups}
We check in this subsection, that every local toric modification
is a composition of a finite number of blowing-ups. 
It will prepare the ground for an inductive proof of Theorem~\ref{T3}.

Let us start from a lemma describing mutual position of vectors spanning a fan.

\begin{Lemma}\label{L1}
If a fan $\cal W$ is spanned by vectors
$\vec\xi_0$, $\vec\xi_1$, \dots , $\vec\xi_n$ and vectors
$\vec\xi_k$, $\vec\xi_l$ ($k+1<l$) form a base of a lattice
then one of vectors~$\vec\xi_i$ is equal to $\vec\xi_k+\vec\xi_l$.
\end{Lemma}

\textbf{Proof.} Suppose that this is not the case. Then $\vec\xi_k+\vec\xi_l$
is inside one of cones~$\sigma_j$ of fan $\cal W$ where $k<j\leq l$ and at least
one of vectors generating $\sigma_j$ is different from $\vec\xi_k$ and $\vec\xi_l$.
We may assume without loss of generality that this is $\vec\xi_j$. Therefore we have the 
following equations with integer coefficients:
$$
\begin{array}{ccl}
 \vec\xi_{j-1} &=& n_1\vec\xi_k + n_2\vec\xi_l \\
 \vec\xi_j     &=& m_1\vec\xi_k + m_2\vec\xi_l \\
 \vec\xi_k + \vec\xi_l  &=& a\vec\xi_{j-1} + b\vec\xi_j \\
 \end{array}
 $$
with $n_1\geq0$, $n_2\geq0$, $n_1+n_2\geq 1$, $m_1>0$, $m_2>0$, $a>0$, $b>0$.
Substituting right-hand sides of two first equations to the third we conclude that
$an_1+bm_1=1$ and $an_2+bm_2=1$, so $n_1=n_2=0$ and we arrive to a contradiction.

\begin{Corollary}\label{W1}
If a fan ${\cal W}_n$ is spanned by vectors
$\vec\xi_0$, $\vec\xi_1$, \dots , $\vec\xi_n$ then for some $i$ ($1<i<n$) we have
$$\vec\xi_i=\vec\xi_{i-1}+\vec\xi_{i+1}.$$
\end{Corollary}

The proof of a corollary uses a simple recurrence. Vectors $\vec\xi_0$ and
$\vec\xi_n$ satisfy assumptions of Lemma~\ref{L1}. Hence among vectors
$\vec\xi_i$ there is one of the form $\vec\xi_0+\vec\xi_n$. Next we apply Lemma~\ref{L1} to a pair 
$\vec\xi_0$, $\vec\xi_i$ or to a pair $\vec\xi_i$, $\vec\xi_n$. Continuing this procedure we 
arrive to such vectors $\vec\xi_{i-1}$, $\vec\xi_{i}$, $\vec\xi_{i+1}$, that
$\vec\xi_i=\vec\xi_{i-1}+\vec\xi_{i+1}$.

It follows from Corollary~\ref{W1} that for every fan ${\cal W}_{n+1}$ consisting of $n+1$ cones 
there exist a fan ${\cal W}_{n}$ consisting of $n$ cones such that one of cones
$\sigma=\sigma[\vec\xi,\vec\nu]\in{\cal W}_{n}$ decomposes into two cones
$\sigma[\vec\xi,\vec\xi+\vec\nu]$, $\sigma[\vec\xi+\vec\nu,\vec\nu]\in{\cal W}_{n+1}$.
We call ${\cal W}_{n+1}$ a subdivision of ${\cal W}_n$.

In the following theorem we compare local toric modification associated 
with a fan and with its subdivision.

\begin{Theorem}\label{T4}
Let $\pi_n:M_n\to\Cn^2$, $\pi_{n+1}:M_{n+1}\to\Cn^2$
be toric modifications associated with a fan ${\cal W}_n$ and its subdivision ${\cal W}_{n+1}$.
Let $\tilde\sigma_i$ be the cone of ${\cal W}_n$ which is divided into two. 
Then $\pi_{n+1}=\pi_n\circ\sigma$ where $\sigma$ is a blowing-up of~$M_n$ at the origin
of the local coordinate system associated with~$\tilde\sigma$.
\end{Theorem}

\textbf{Proof.} Let 
$\tilde\sigma=\sigma[\vec\xi,\vec\nu]$, 
$\sigma'=\sigma[\vec\xi,\vec\xi+\vec\nu]$, 
$\sigma''=\sigma[\vec\xi+\vec\nu,\vec\nu]$. 
In all charts associated with cones 
different from $\tilde\sigma$, $\sigma'$, $\sigma''$
mappings $\pi_n$ and $\pi_{n+1}$ are given by identical formulas, hence in these charts
$\sigma$ is an identity. If we examine $\varphi_{\sigma'}$,
$\varphi_{\sigma''}$ and $\varphi_{\tilde\sigma}$ then it is easy to check
that $\varphi_{\sigma'}=\varphi_{\tilde\sigma}\circ\sigma$ where $\sigma(u,v)=(uv,v)$ and
$\varphi_{\sigma''}=\varphi_{\tilde\sigma}\circ\sigma$ where $\sigma(u,v)=(u,uv)$.
Hence $\sigma$ is a blowing-up of $M_n$ at a point $(0,0)$ of a chart associated with $\tilde\sigma$.

\subsection{Orders of proper preimages in centers of blowing-ups}

\begin{Lemma}\label{Lr}
Let $f=f(x,y)$ be a convenient power series with a Newton diagram $\Delta$.
Let $\pi:M\to\Cn^2$ be a local toric modification with a fan $\cal W$ and let $\sigma$ be a cone of~$\cal W$ 
spanned by vectors $\vec\xi$, $\vec\nu$. Then the order of the proper preimage of~$f$ at a point~$(0,0)$
of a chart associated with~$\sigma$ is equal to
$$
l(\Delta,\vec\xi+\vec\nu) - l(\Delta,\vec\xi) - l(\Delta,\vec\nu) .
$$
\end{Lemma}

\textbf{Proof.}
For
$$ f(x,y) = \sum_{ij}a_{ij}x^iy^j $$
we have
$$ (f\circ\varphi_\sigma)(u,v) = \sum_{ij}a_{ij}u^{\scalar{(i,j)}{\vec\xi}}v^{\scalar{(i,j)}{\vec\nu}}
$$
We may exclude factors $u^{l(\Delta,\vec\xi)}$ and $v^{l(\Delta,\vec\nu)}$ from above sum getting
\begin{equation}\label{L4.1}
f(x,y) = u^{l(\Delta,\vec\xi)}v^{l(\Delta,\vec\nu)}\tilde f(u,v)
\end{equation}
where $\tilde f(u,v)$ is a proper preimage of $f$.

The order of $\tilde f$ at $(0,0)$ is equal to the order of substitution
of generic curve $\gamma:t\to(c_1t,c_2t)$ to~$\tilde f$. In coordinates $(x,y)$
the curve $\gamma$ has an equation
$(x,y)=(d_1t^{\xi_1+\nu_1},d_2t^{\xi_2+\nu_2})$ and for generic $d_1$, $d_2$ we have
$\ord f(d_1t^{\xi_1+\nu_1},d_2t^{\xi_2+\nu_2})= l(\Delta,\vec\xi+\vec\nu)$ (see~\cite{Kh}).
Thus substituting the parameterization of $\gamma$ to (\ref{L4.1}) we get
$$
l(\Delta,\vec\xi+\vec\nu) =
\ord\bigl(t^{l(\Delta,\vec\xi)}t^{l(\Delta,\vec\nu)}(\tilde f\circ\gamma)(t)\bigr) =
l(\Delta,\vec\xi) + l(\Delta,\vec\nu) + \ord\tilde f.
$$

\subsection{Proof of Theorem~\ref{T3}}

\vspace{1ex}
\noindent
\textbf{Proof of~(\ref{T3.2}).}
Let ${\cal W}_n$, ${\cal W}_{n+1}$ be fans from Theorem~\ref{T4}.
Let $f$, $g$ be convenient power series such that
$\tilde\Delta_f=\tilde\Delta_g=\tilde\Delta$ where thickened Newton diagrams are relative to
${\cal W}_{n+1}$. Keeping notation from Theorem~\ref{T4} we shall check that from formula
\begin{equation}\label{D1}
\delta_0(f)- \sum_{p\in\pi_n^{-1}(0)} \delta_p(\tilde f_n) \quad=\quad
\delta_0(g)- \sum_{p\in\pi_n^{-1}(0)} \delta_p(\tilde g_n)
\end{equation}
follows
\begin{equation}\label{D2}
\delta_0(f)- \sum_{p\in\pi_{n+1}^{-1}(0)} \delta_p(\tilde f_{n+1}) \quad=\quad
\delta_0(g)- \sum_{p\in\pi_{n+1}^{-1}(0)} \delta_p(\tilde g_{n+1})
\end{equation}

We use subscripts for $\tilde f_n$, $\tilde g_n$, $\tilde f_{n+1}$, $\tilde g_{n+1}$
to distinguish proper preimages on manifolds $M_n$ and $M_{n+1}$.
By Theorem~\ref{T4} $\pi_{n+1}=\pi_n\circ\sigma$ where $\sigma$ is
a blowing-up of $M_n$ at the center $q$ which is the origin of the chart associated with 
$\tilde\sigma=\sigma[\vec\xi,\vec\nu]$. 
By Lemma~\ref{Lr} orders of $\tilde f_n$ and $\tilde g_n$ at $q$
are identical and equal to
$l(\tilde\Delta,\vec\xi+\vec\nu) - l(\tilde\Delta,\vec\xi) - l(\tilde\Delta,\vec\nu)$.
Therefore from Noether formula~(\ref{T2.2}) we get
\begin{equation}\label{D3}
\delta_q(\tilde f_n) - \sum_{p\in\sigma^{-1}(q)} \delta_p(\tilde f_{n+1}) \quad=\quad
\delta_q(\tilde g_n) - \sum_{p\in\sigma^{-1}(q)} \delta_p(\tilde g_{n+1})
\end{equation}
Adding (\ref{D1}) and (\ref{D3}) we get (\ref{D2}).

An inductive argument with respect to the number of cones in the fan leads to conclusion
that formula (\ref{D1}) holds for every fan $\cal W$
and all convenient power series $f$, $g$ such that
$\tilde\Delta_f=\tilde\Delta_g=\tilde\Delta$.
Taking as $g$ a nondegenerate convergent power series with a Newton diagram $\tilde\Delta$
and using Theorem~\ref{T1}
we get $\delta_0(g) = \delta(\tilde\Delta)$ and
$\delta_p(\tilde g_n)=0$ for every $p\in\pi_n^{-1}(0)$.
Hence the right-hand side of~(\ref{D1}) is equal to~$\delta(\tilde\Delta)$ which  proves~(\ref{T3.2}).

\vspace{2ex}
\noindent
\textbf{Proof of~(\ref{T3.3}).}
It is enough to use a formula $\delta_p(f)=(\mu_p(f)+r_p(f)-1)/2$.
After substituting to~(\ref{T3.2}) and multiplying by 2 we get
$$
\mu_0(f)+r_0(f)-1 = \mu(\tilde\Delta)+r(\tilde\Delta)-1+ \sum_{p\in\pi^{-1}(0)\cap\{\tilde f=0\}}
(\mu_p(\tilde f )+r_p(\tilde f)-1)
$$
Then we apply an equality on number of branches
$r_0(f) =\sum_{p\in\pi^{-1}(0)} r_p(\tilde f)$
and we get~(\ref{T3.3}).

\vspace{2ex}
\noindent
\textbf{Proof of~(\ref{T3.1}).} 
We only outline the proof briefly because it is analogous to the one of~(\ref{T3.3}).
First we show by induction on the number of cones that
\begin{equation}\label{D4}
(f_1,g_1)_0 - \sum_{p\in\pi^{-1}(0)} (\tilde f_1,\tilde g_1)_p \quad=\quad
(f_2,g_2)_0 - \sum_{p\in\pi^{-1}(0)} (\tilde f_2,\tilde g_2)_p
\end{equation}
where $f_i$, $g_i$ ($i=1,2$) are such convergent power series that
$\tilde\Delta_{f_1}=\tilde\Delta_{f_2}$ and $\tilde\Delta_{g_1}=\tilde\Delta_{g_2}$.
In the inductive proof we use a formula~(\ref{T2.1})
for the intersection number after blow-up.

Then we take as $f_2$ and $g_2$ nondegenerate power series
with Newton diagrams $\Delta_{f_2}=\tilde\Delta_{f_1}$ and
$\Delta_{g_2}=\tilde\Delta_{g_1}$. By Theorem~\ref{T1} 
the right-hand side of~(\ref{D4}) is equal to
$[\Delta_{f_2},\Delta_{g_2}]$ which proves~(\ref{T3.1}).

\end{document}